\documentclass[12pt]{article}
\usepackage{amsmath}
\usepackage{amsfonts}
\usepackage{mathrsfs}
\usepackage{amssymb}
\usepackage[usenames]{color}
\def\dis{\displaystyle}

\def\ssc{\scriptscriptstyle}
\def\OP#1{\raisebox{-5pt}{$\stackrel{\dis\oplus}{\ssc #1}$}}

\def\cl{\centerline}

\def\a{\alpha}

\def\b{\beta}
\def\vs{\vspace*}

\def\Z{\mathbb{Z}}

\def\C{\mathbb{C}}

\def\QED{\hfill$\Box$}
\def\ni{\noindent}

\numberwithin{equation}{section}
\newtheorem{theo}{Theorem}[section]
\newtheorem{defi}[theo]{Definition}

\newtheorem{lemm}[theo]{Lemma}
\newtheorem{rema}[theo]{Remark}

\begin{document}
\begin{center}
{\bf\large Infinite rank Schr$\ddot{\rm o}$dinger-Virasoro type Lie conformal algebras}
\footnote{
$^{\,\dag}$Corresponding author: chgxia@cumt.edu.cn (C.~Xia).}
\end{center}

\cl{Guangzhe Fan$^{\,*}$, \ Yucai Su$^{\,*}$, \ Chunguang Xia$^{\,\dag}$}
\vspace{10pt}

\cl{\small $^{\,*}$Department of Mathematics, Tongji University, Shanghai 200092, China}
\vspace{10pt}

\cl{\small $^{\,\dag}$Department of Mathematics, China University of Mining and Technology, Xuzhou 221116, China}
\vspace{10pt}

\cl{\small Email: yzfanguangzhe@126.com, ycsu@tongji.edu.cn, chgxia@cumt.edu.cn}
\vs{10pt}

\small\footnotesize
\noindent{{\bf Abstract.}
Motivated by the structure of certain modules over the loop Virasoro Lie conformal algebra and
the Lie structures of Schr$\ddot{\rm o}$dinger-Virasoro algebras, we construct a class of infinite rank Lie conformal algebras $CSV(a,b)$,
where $a,\,b$ are complex numbers.
The conformal derivations of $CSV(a,b)$ are uniformly determined.
The rank one conformal modules and $\Z$-graded free intermediate series modules over $CSV(a,b)$ are classified.
Corresponding results of the conformal subalgebra $CHV(a,b)$ of $CSV(a,b)$ are also presented.
\vs{5pt}

\ni{\bf Key words:} Lie conformal algebra, conformal derivation, conformal module

\ni{\it Mathematics Subject Classification (2010):} 17B05; 17B10; 17B40; 17B65; 17B68.}

\section{Introduction}
The notion of Lie conformal algebra, introduced by Kac \cite{K3}, encodes the singular parts of the operator
product expansion in conformal field theory.
Conformal modules are basic tools for the construction of free field realization of infinite dimensional Lie (super)algebras in conformal field theory.
In the language of the $\lambda$-bracket, a Lie conformal algebra $R$ is a $\C[\partial]$-module endowed with a
$\C$-linear map $R\otimes R\rightarrow R[\lambda]$, where $R[\lambda]=\C[\lambda]\otimes R$, satisfying the following axioms ($a,\,b,\,c\in R$):
\begin{equation*}
\aligned
\mbox{(conformal sesquilinearity)}&~~~~[\partial a\,{}_\lambda \,b]=-\lambda[a\,{}_\lambda\, b],\ \ \ \
[a\,{}_\lambda \,\partial b]=(\partial+\lambda)[a\,{}_\lambda\, b],\\
\mbox{(skew-symmetry)}&~~~~[a\, {}_\lambda\, b]=-[b\,{}_{-\lambda-\partial}\,a],\\
\mbox{(Jacobi identity)}&~~~~[a\,{}_\lambda\,[b\,{}_\mu\, c]]=[[a\,{}_\lambda\, b]\,{}_{\lambda+\mu}\, c]+[b\,{}_\mu\,[a\,{}_\lambda \,c]].
\endaligned
\end{equation*}
In recent years, the structure and representation theory of Lie conformal algebras have been extensively studied.
Finite simple Lie conformal algebras were classified in \cite{DK} and their cohomology theory
and representation theory were further developed in \cite{BKV,CK}.
The simplest but rather important finite Lie conformal algebra is the Virasoro conformal algebra, whose irreducible
conformal modules were classified in \cite{CK}.
More recently, some finite Lie conformal algebras related to the Virasoro conformal algebra have been constructed and studied,
such as Schr$\ddot{\rm o}$dinger-Virasoro type Lie conformal algebras \cite{H,SYua,WXX} and $W(a,b)$ Lie conformal algebras \cite{XY}.
Infinite rank Lie conformal algebras are also important ingredients of Lie conformal algebras, and have
also attracted many authors' attention, such as general Lie conformal algebras \cite{BKL1,BKL2,S,SYue},
and some loop Lie conformal algebras \cite{CHSX,FSW,FWY,WCY}.

Motivated by the structure of certain modules over the loop Virasoro Lie conformal algebra \cite{WCY} and
the Lie structures of Schr$\ddot{\rm o}$dinger-Virasoro algebras \cite{He,RU},
in this paper, we construct a class of infinite rank Lie conformal algebras $CSV(a,b)$, where $a,b\in\C$.
In view of our construction (see Section~2 for details), we call these conformal algebras
{\it infinite rank Schr$\ddot{o}$dinger-Virasoro type Lie conformal algebras}.
Precisely speaking, $CSV(a,b)$ has a $\C[\partial]$-basis $\{L_i,M_i,Y_i\,|\,i\in\Z\}$
satisfying the following nonvanishing $\lambda$-brackets (see Theorem~\ref{theorem-construction}):
\begin{eqnarray}
\label{brackets-1}[L_i\, {}_\lambda \, L_j] &\!\!=\!\!& (\partial+2\lambda) L_{i+j},\\
\label{brackets-2}[L_i\, {}_\lambda \, M_j] &\!\!=\!\!& (\partial+a\lambda+b) M_{i+j}, \\
\label{brackets-3}[L_i\, {}_\lambda \, Y_j] &\!\!=\!\!& \mbox{$(\partial+(\frac{a}{2}+1)\lambda+\frac{b}{2}) Y_{i+j}$}, \\
\label{brackets-4}[Y_i\, {}_\lambda \, Y_j] &\!\!=\!\!& (\partial+2\lambda)M_{i+j}.
\end{eqnarray}
Note that the special case $CSV(1,0)$ is exactly the loop Schr$\ddot{\rm o}$dinger-Virasoro Lie conformal algebra recently studied in \cite{CHSX}.
In addition, $CSV(a,b)$ contains many important conformal subalgebras.
For example,
\begin{itemize}
  \item $CW={\oplus}_{i\in\Z}\C[\partial]L_i$ is in fact the loop
  Virasoro Lie conformal algebra introduced in \cite{WCY}.

  \item $CHV(a,b)=({\oplus}_{i\in\Z}\C[\partial]L_i)\oplus({\oplus}_{i\in\Z}\C[\partial]M_i)$
  is an infinite rank conformal subalgebra, which we call an {\it infinite rank Heisenberg-Virasoro type Lie conformal algebra}.
  The special case $CHV(a,0)$ is in fact the loop Lie conformal algebra $R(1-a)$ studied in \cite{FWY},
  and the more special case $CHV(1,0)$ is exactly the loop Heisenberg-Virasoro Lie conformal algebra firstly constructed in \cite{FSW}.

  \item $SV(a,b)=\C[\partial]L_0\oplus\C[\partial]M_0\oplus\C[\partial]Y_0$ is a {\it finite
  Schr$\ddot{o}$dinger-Virasoro type conformal subalgebra}, which was studied in \cite{H}.
  The special cases $SV(1,0)$ and $SV(-2,0)$ were firstly considered in \cite{SYua} and \cite{WXX}, respectively.
  Furthermore, $SV(a,b)$ contains a Virasoro conformal subalgebra ${\rm CVir}=\C[\partial]L_{0}$, and
  a {\it finite Heisenberg-Virasoro type conformal subalgebra} $HV(a,b)=\C[\partial]L_0\oplus\C[\partial]M_0$. 
  The special case $HV(a,0)$ was studied in \cite{XY}, and the more special case $HV(1,0)$ is exactly the Heisenberg-Virasoro conformal algebra
  introduced in \cite{SYua}.
\end{itemize}
Due to the above facts, it seems to be interesting for us to study  the structure and representation theory of the parent algebra $CSV(a,b)$ in a uniform way.
 This is the main motivation of this paper.

We will study the conformal derivations, rank one conformal modules and
$\Z$-graded free intermediate series modules for $CSV(a,b)$ and its conformal subalgebra $CHV(a,b)$ for all $a,b\in\C$.
Comparing our results with those for $CW$ in \cite{WCY}, we find some nontrivial conclusions
(see Table~1 and Table~2).
$$
\begin{tabular}{c|c|c}
\hline
 $CSV(a,b)$  & Values of $a$, $b$ & Reference \\
\hline
Non-inner conformal derivations  & $a=1$ & Theorem~\ref{main-results-1} \\
\hline
Nontrivial extentions from $CW$-modules
& $a=b=0$ &
Theorems~\ref{main-results-2} and \ref{main-results-3}
\\
\hline
\end{tabular}
$$
\vs{-25pt}
\begin{center}
\mbox{Table 1: Necessary conditions for existence of nontrivial conclusions on $CSV(a,b)$}
\end{center}
\vs{-20pt}
$$
\begin{tabular}{c|c|c}
\hline
 $CHV(a,b)$  & Values of $a$, $b$ & Reference \\
\hline
Non-inner conformal derivations  & $a=1$ & Theorem~\ref{main-results-1-CHV} \\
\hline
Nontrivial extentions from $CW$-modules
& $a=1,\,b=0$ &
Theorems~\ref{main-results-2-CHV} and \ref{main-results-3-CHV}
\\
\hline
\end{tabular}
$$
\vs{-20pt}
\begin{center}
\mbox{Table 2: Necessary conditions for existence of nontrivial conclusions on $CHV(a,b)$}
\end{center}

This paper is arranged as follows.
In Section~2, we give the detail construction of the Lie conformal algebra $CSV(a,b)$.
Then, in Section~3, we uniformly determine the conformal derivations of $CSV(a,b)$.
Next, in Sections~4 and 5, we classify the rank one conformal modules and
$\Z$-graded free intermediate series modules over $CSV(a,b)$, respectively.
Finally, we present the corresponding results of the conformal subalgebra $CHV(a,b)$ of $CSV(a,b)$ in Section~6.

\section{Lie conformal algebras $CSV(a,b)$}

All definitions related to Lie conformal algebras used in this paper are collected from \cite{DK,K1,K3}.
In this section, we give the construction of the Lie conformal algebra $CSV(a,b)$ for $a,b\in\C$.
First, let us recall some definitions on conformal modules.
\begin{defi}\label{def-module}\rm
A conformal module $M$ over a Lie conformal algebra $R$ is a $\C[\partial]$-module endowed with a $\lambda$-action
$R\otimes M\rightarrow M[\lambda]$ such that ($a,\,b\in R$, $v\in V$)
\begin{equation*}
(\partial a)\,{}_\lambda\, v=-\lambda a\,{}_\lambda\, v,\ \ \ \ \ a{}\,{}_\lambda\, (\partial v)=(\partial+\lambda)a\,{}_\lambda\, v,\ \ \ \
a\,{}_\lambda\, (b{}\,_\mu\, v)-b\,{}_\mu\,(a\,{}_\lambda\, v)=[a\,{}_\lambda\, b]\,{}_{\lambda+\mu}\, v.
\end{equation*}
\end{defi}
\begin{defi}\rm
A Lie conformal algebra $R$ is {\it $\Z$-graded} if $R=\oplus_{i\in \Z}R_i$, where each $R_i$ is a $\C[\partial]$-submodule
and $[R_i\,{}_\lambda\,R_j]\subset R_{i+j}[\lambda]$ for $i,j\in \Z$.
Similarly, a conformal module $V$ over $R$ is  {\it $\Z$-graded} if $V=\oplus_{i\in \Z}V_i$, where each $V_i$ is a $\C[\partial]$-submodule and
$R_i\,{}_\lambda\, V_j\subset V_{i+j}[\lambda]$ for $i,j\in \Z$. Furthermore, if each $V_i$ is freely generated by one element $v_i\in V_i$
over $\C[\partial]$, we call $V$ a {\it $\Z$-graded free intermediate series module}.
\end{defi}

The start point of our construction is the loop Virasoro Lie conformal algebra $CW={\oplus}_{i\in\Z}\C[\partial]L_i$
with $\lambda$-brackets \eqref{brackets-1}. According to \cite{WCY},
one class of $\Z$-graded free intermediate series modules over $CW$ is $V_{a,b}$, $a,b\in\C$, which has $\C[\partial]$-basis $\{v_j\,|\,j\in\Z\}$ and the $\lambda$-actions are given by
\begin{equation}\label{CW-actions}
L_i\,{}_\lambda\, v_j=(\partial+a\lambda+b)v_{i+j}\ \ \ \mbox{for}\ \ \ i,j\in\Z.
\end{equation}
The original Schr$\ddot{\rm o}$dinger-Virasoro algebra $\mathfrak{sv}=\text{span}_\C\{L_m,\,M_n,\,Y_p\,|\,m,n\in\Z,\,p\in \frac{1}{2}+\Z\}$ was introduced by Henkel in \cite{He}
in the context of nonequilibrium statistical physics, while the twisted Schr$\ddot{\rm o}$dinger-Virasoro algebra
$\mathfrak{tsv}=\text{span}_\C\{L_m,\,M_n,\,Y_p\,|\,m,n\in\Z,\,p\in \Z\}$ was introduced by Roger and Unterberger in \cite{RU}.
The Lie algebras $\mathfrak{sv}$ and $\mathfrak{tsv}$ have the same Lie structures, which are given by (others vanishing)
\begin{eqnarray}
\label{sv-brackets-1}[L_m, L_n] &\!\!=\!\!& (n-m) L_{m+n},\\
\label{sv-brackets-2}[L_m, M_n] &\!\!=\!\!& n M_{m+n}, \\
\label{sv-brackets-3}[L_m, Y_p] &\!\!=\!\!& \mbox{$(p-\frac{m}{2}) Y_{m+p}$}, \\
\label{sv-brackets-4}[Y_p, Y_{p'}] &=& (p'-p)M_{p+p'}.
\end{eqnarray}
Note that ${\rm Vir}=\text{span}_\C\{L_m\,|\,m\in\Z\}$ is a Virasoro subalgebra of $\mathfrak{sv}$ (resp.~$\mathfrak{tsv}$).
The formulas \eqref{sv-brackets-2} and \eqref{sv-brackets-3} say that the spaces spanned respectively by $M_n$ and $Y_p$
are modules over ${\rm Vir}$.
For $a,a',b,b'\in\C$, motivated by \eqref{CW-actions} and \eqref{sv-brackets-1}--\eqref{sv-brackets-4},
it is nature to define a  free $\C[\partial]$-module
$M(a,a',b,b')$, which has $\C[\partial]$-basis $\{L_i,M_i,Y_i\,|\,i\in\Z\}$ and satisfies the following nonvanishing $\lambda$-brackets:
\begin{eqnarray}
\nonumber[L_i\, {}_\lambda \, L_j]&\!\!=\!\!&(\partial+2\lambda) L_{i+j},\\
\nonumber[L_i\, {}_\lambda \, M_j]&\!\!=\!\!&(\partial+a\lambda+b) M_{i+j}, \\
\nonumber[L_i\, {}_\lambda \, Y_j]&\!\!=\!\!&(\partial+a'\lambda+b') Y_{i+j},\\
\nonumber[Y_i\, {}_\lambda \, Y_j]&\!\!=\!\!&(\partial+2\lambda)M_{i+j}.
\end{eqnarray}

\begin{theo}\label{theorem-construction}
$M(a,a',b,b')$ is a Lie conformal algebra if and only if $a'=\frac{a}{2}+1$ and $b'=\frac{b}{2}$.
\end{theo}
\ni{\it Proof.}\ \
Suppose $M(a,a',b,b')$ is a Lie conformal algebra. By the Jacobi identity, we have
$$[L_{i}\,{}_\lambda\,[Y_{j}\,{}_\mu\, Y_{k}]]=[[L_{i}\,{}_\lambda\, Y_{j}]\,{}_{\lambda+\mu}\, Y_{k}]+[Y_{j}\,{}_\mu\,[L_{i}\,{}_\lambda \,Y_{k}]].$$
A direct computation shows that
\begin{equation*}
(\partial+\lambda+2\mu)(\partial+a\lambda+b)=(\partial+2\lambda+2\mu)((a'-1)\lambda-\mu+b')+(\partial+2\mu)(\partial+a'\lambda+\mu+b').
\end{equation*}
Comparing the coefficients of $\partial\lambda$ and $\partial$,
one can deduce that
$a'=\frac{a}{2}+1$ and $b'=\frac{b}{2}$.
On the other hand, if $a'=\frac{a}{2}+1$ and $b'=\frac{b}{2}$,
then $M(a,a',b,b')$ is indeed a Lie conformal algebra (we leave the verification details to the reader).
\QED

\begin{rema}\rm
(1) The above theorem gives the necessary and sufficient conditions for $M(a,a',b,b')$ to be
a Lie conformal algebra.
In view of the above construction, we call $M(a,a',b,b')$ an
{\it infinite rank Schr$\ddot{o}$dinger-Virasoro type Lie conformal algebra},
and redenote it by $CSV(a,b)$. 

(2) Note that the Lie subalgebra with $\C$-basis $\{L_m,\,M_n\,|\,m,n\in\Z\}$ of $\mathfrak{sv}$ (or $\mathfrak{tsv}$)
is in fact the twisted Heisenberg-Virasoro algebra introduced by Arbarello et al. in \cite{ADKP}.
Naturally, we call the Lie conformal subalgebra with $\C[\partial]$-basis $\{L_i,\,M_i\,|\,i\in\Z\}$
of $CSV(a,b)$ an {\it infinite rank Heisenberg-Virasoro type Lie conformal algebra}, and denote it by $CHV(a,b)$.

(3) Note that $CSV(a,b)$ is a $\Z$-graded Lie conformal algebra in the sense
$CSV(a,b)=\oplus_{i\in\Z} (CSV(a,b))_i$, where
$(CSV(a,b))_i=\C[\partial]L_i \oplus \C[\partial] M_i \oplus \C[\partial] Y_i$.
\end{rema}

\section{Conformal derivations of $CSV(a,b)$}

\begin{defi}\rm
Let $R$ be a Lie conformal algebra. A linear map $D_\lambda:R\rightarrow R[\lambda]$ is called a conformal derivation if ($a,\,b\in R$)
$$
D_\lambda(\partial a)=(\partial+\lambda)D_\lambda(a),\ \ \ \
D_\lambda([a\,{}_\mu \,b])=[(D_\lambda a)\,{}_{\lambda+\mu} \,b]+[a\,{}_\mu \,(D_\lambda b)].
$$
We often write $D$ instead of $D_\lambda$ for simplicity.
By the Jacobi identity, for any $a\in R$, the map ${\rm ad}_a$, defined by $({\rm ad}_a)_\lambda b= [a\, {}_\lambda\, b]$ for $b\in R$, is a conformal derivation of $R$. All derivations of this kind are called {\it inner}.
\end{defi}

Denote by ${\rm CDer\,}(CSV(a,b))$
and ${\rm CInn\,}(CSV(a,b))$ the vector spaces of all conformal derivations and inner conformal derivations of $CSV(a,b)$, respectively.
From \cite{WCY}, all conformal derivations of the loop  Virasoro Lie conformal algebra $CW$ are inner.
We focus our interest on whether there are non-inner conformal derivations of $CSV(a,b)$. Let us consider $CSV(1,b)$. Denote
$$
\C^\infty=\{\vec{a}=(a_c)_{c\in\Z}\,|\,a_c\in\C\mbox{ and }a_c=0\ \mbox{for\ all\ but\ finitely\ many}\ c\mbox{'}s\}.
$$
For each $\vec{a}\in \C^\infty$, we define $D_{\vec{a}}$ by
\begin{equation}\label{non-inner}
D_{\vec{a}}\, {}_\lambda\, (L_i)=\sum_{c\in\Z} a_c M_{i+c}, \ \ \ D_{\vec{a}}\, {}_\lambda\,(M_i)=D_{\vec{a}}\, {}_\lambda\,(Y_i)=0 \ \ \ \mbox{for}\ \ \  i\in\Z.
\end{equation}
Then $D_{\vec{a}}\in {\rm CDer\,}(CSV(1,b))$. Furthermore, if $D_{\vec{a}}\in {\rm CInn\,}(CSV(1,b))$, then one can easily show that $D_{\vec{a}}=0$.
For simplicity, we still use $\C^\infty$ to stand for the space of such conformal derivations.
We will show that there are no other non-inner conformal derivations of $CSV(a,b)$ except non-zero $D_{\vec{a}}$'s defined in \eqref{non-inner}.

\begin{theo}\label{main-results-1}
 We have ${\rm CDer\,}(CSV(a,b))={\rm CInn\,}(CSV(a,b))\oplus \delta_{a,1}\C^\infty$.
\end{theo}


First, we prove the following lemma.
For $c\in\Z$, denote by ${({\rm CDer\,}(CSV(a,b))}^c$ the space of conformal derivations of degree $c$, i.e.,
$$
{({\rm CDer\,}(CSV(a,b)))}^c=\{D\in {\rm CDer\,}(CSV(a,b))\,|\, D_\lambda((CSV(a,b))_i)\subset (CSV(a,b))_{i+c}[\lambda]\}.
$$

\begin{lemm}\label{graded-derivations}
Any $D^c\in {({\rm CDer\,}(CSV(a,b)))}^c$ can be written as $D^c={\rm ad}_x+\delta_{a,1}D_{\vec{a}}$ for some $x\in (CSV(a,b))_c$ and
$\vec{a}=(\ldots, 0, \ldots, 0, q, 0, \ldots, 0, \ldots)\in\C^\infty$, where $q\in\C$ appears in the $c$-th position.
\end{lemm}

\ni{\it Proof.}\ \
Fix $c\in\Z$. Take $D^c\in {({\rm CDer\,}(CSV(a,b)))}^c$.
We only need to show that $D^c$ minus some suitable derivations (inner derivations or $D_{\vec{a}}$) equal to zero derivation.
Assume that
\begin{eqnarray*}
D^c_\lambda(L_i)&\!\!=\!\!&f_{1,i}(\partial,\lambda)L_{i+c}+f_{2,i}(\partial,\lambda)M_{i+c}+f_{3,i}(\partial,\lambda)Y_{i+c},\\
D^c_\lambda(M_i)&\!\!=\!\!&g_{1,i}(\partial,\lambda)L_{i+c}+g_{2,i}(\partial,\lambda)M_{i+c}+g_{3,i}(\partial,\lambda)Y_{i+c},\\
D^c_\lambda(Y_i)&\!\!=\!\!&h_{1,i}(\partial,\lambda)L_{i+c}+h_{2,i}(\partial,\lambda)M_{i+c}+h_{3,i}(\partial,\lambda)Y_{i+c}.
\end{eqnarray*}
Applying $D^c_\lambda$ to $[L_0\ {}_\mu \ L_i]=(\partial+2\mu)L_i$ and equating coefficients, we obtain
\begin{eqnarray}
\nonumber
(\partial+\lambda+2\mu)f_{1,i}(\partial,\lambda)&\!\!=\!\!&(\partial+2\lambda+2\mu)f_{1,0}(-\lambda-\mu,\lambda)+(\partial+2\mu)f_{1,i}(\partial+\mu,\lambda),\\
\nonumber
(\partial+\lambda+2\mu)f_{2,i}(\partial,\lambda)&\!\!=\!\!&((a-1)\partial+a(\lambda+\mu)-b)f_{2,0}(-\lambda-\mu,\lambda)+(\partial+a\mu+b)f_{2,i}(\partial+\mu,\lambda),\\
\nonumber
(\partial+\lambda+2\mu)f_{3,i}(\partial,\lambda)&\!\!=\!\!&\mbox{$(\frac{a}{2}\partial+(\frac{a}{2}+1)(\lambda+\mu)-\frac{b}{2})f_{3,0}(-\lambda-\mu,\lambda)+
(\partial+(\frac{a}{2}+1)\mu+\frac{b}{2})f_{3,i}(\partial+\mu,\lambda)$}.
\end{eqnarray}
In particular, putting $\mu=0$, we have
\begin{eqnarray}
\label{404}
\lambda f_{1,i}(\partial,\lambda)&\!\!=\!\!&(\partial+2\lambda)f_{1,0}(-\lambda,\lambda),\\
\label{405}
(\lambda-b)f_{2,i}(\partial,\lambda)&\!\!=\!\!&((a-1)\partial+a\lambda-b)f_{2,0}(-\lambda,\lambda),\\
\label{406}
\mbox{$(\lambda-\frac{b}{2})f_{3,i}(\partial,\lambda)$}&\!\!=\!\!&\mbox{$(\frac{a}{2}\partial+(\frac{a}{2}+1)\lambda-\frac{b}{2})f_{3,0}(-\lambda,\lambda).$}
\end{eqnarray}
By \eqref{404}, $\lambda$ is a factor of  $f_{1,0}(-\lambda,\lambda)$ in the unique factorization ring $\C[\partial,\lambda]$.
Setting $f_1(\lambda)=\frac{f_{1,0}(-\lambda,\lambda)}{\lambda}$ 
and replacing $D^c$ by $D^c-{\rm ad}_{{f_1(-\partial)}L_c}$,
we may assume that $D^c_\lambda(L_i)=f_{2,i}(\partial,\lambda)M_{i+c}+f_{3,i}(\partial,\lambda)Y_{i+c}$.
We divide the remaining proof into the following three cases.

\vskip2pt
{\bf Case 1:} $a=1$.
\vskip2pt

In this case, \eqref{405} and \eqref{406} become
\begin{eqnarray}
\label{407}
(\lambda-b) f_{2,i}(\partial,\lambda)&\!\!=\!\!&(\lambda-b) f_{2,0}(-\lambda,\lambda),\\
\label{408}
\mbox{$(\lambda-\frac{b}{2}) f_{3,i}(\partial,\lambda)$}&\!\!=\!\!&\mbox{$(\frac{1}{2}\partial+\frac{3}{2}\lambda-\frac{b}{2})f_{3,0}(-\lambda,\lambda).$}
\end{eqnarray}
By \eqref{408}, $\lambda-\frac{b}{2}$ is a factor of  $f_{3,0}(-\lambda,\lambda)$ in the unique factorization ring $\C[\partial,\lambda]$.
Replacing $D^c$ by $D^c-{\rm ad}_{f_3(-\partial)Y_c}$ with $f_3(\lambda)=\frac{f_{3,0}(-\lambda,\lambda)}{\lambda-\frac{b}{2}}$,
we may assume that $D^c_\lambda(L_i)=f_{2,i}(\partial,\lambda)M_{i+c}$.
Then, by \eqref{407}, there exists a polynomial $f_2(\lambda)\in\C[\lambda]$ such that $f_{2,i}(\partial,\lambda)=f_2(\lambda)$ for all $i\in\Z$.
By polynomial long division, we may write $f_2(\lambda)=(\lambda-b) f'_2(\lambda)+q$, where $f'_2(\lambda)\in\C[\lambda]$ and $q\in\C$.
For this $q$, we define
$$
\vec{a}=(\ldots, 0, \ldots, 0, q, 0, \ldots, 0, \ldots)\in\C^\infty,
$$
where $q$ appears in the $c$-th position. By \eqref{non-inner}, $D_{\vec{a}}\in\C^\infty$.
Replacing $D^c$ by $D^c-{\rm ad}_{{f'_2(-\partial)}M_c}-D_{\vec{a}}$, we may assume that $D^c_\lambda(L_i)=0$.
Applying $D^c_\lambda$ to $[L_0\ {}_\mu \ M_i]=(\partial+\mu+b)M_i$ and equating coefficients, we obtain
\begin{eqnarray}
\nonumber
(\partial+\lambda+\mu+b)g_{1,i}(\partial,\lambda)&\!\!=\!\!&(\partial+2\mu)g_{1,i}(\partial+\mu,\lambda),\\
\nonumber
(\partial+\lambda+\mu+b)g_{2,i}(\partial,\lambda)&\!\!=\!\!&(\partial+\mu+b)g_{2,i}(\partial+\mu,\lambda),\\
\nonumber
(\partial+\lambda+\mu+b)g_{3,i}(\partial,\lambda)&\!\!=\!\!&\mbox{$(\partial+\frac{3}{2}\mu+\frac{b}{2})g_{3,i}(\partial+\mu,\lambda)$}.
\end{eqnarray}
Putting $\mu=0$ in the above equalities, one can easily deduce that
$g_{1,i}(\partial,\lambda)=g_{2,i}(\partial,\lambda)=g_{3,i}(\partial,\lambda)=0$.
Therefore, $D^c_\lambda(M_i)=0$.
Similarly, one can deduce $D^c_\lambda(Y_i)=0$ by $[L_0\ {}_\mu \ Y_i]=(\partial+\frac{3}{2}\mu+\frac{b}{2})Y_i$.

\vskip2pt
{\bf Case 2:} $a=0$.
\vskip2pt

In this case, \eqref{405} and \eqref{406} become
\begin{eqnarray}
\label{415}
(\lambda-b) f_{2,i}(\partial,\lambda)&\!\!=\!\!&(-\partial-b) f_{2,0}(-\lambda,\lambda),\\
\label{416}
\mbox{$(\lambda-\frac{b}{2}) f_{3,i}(\partial,\lambda)$}&\!\!=\!\!&\mbox{$(\lambda-\frac{b}{2})f_{3,0}(-\lambda,\lambda).$}
\end{eqnarray}
As in Case 1, by \eqref{415} and \eqref{416}, we can set $g_2(\lambda)=\frac{f_{2,0}(-\lambda,\lambda)}{\lambda-b}$
and write $f_{3,i}(\partial,\lambda)=(\lambda-\frac{b}{2})g_3(\lambda)+p$ for all $i\in\Z$, where $g_3(\lambda)\in \C[\lambda]$ and $p\in\C$.
Replacing $D^c$ by $D^c-{\rm ad}_{{g_2(-\partial)}M_c+g_3(-\partial)Y_c}$, we may assume that $D^c_\lambda(L_i)=p Y_{i+c}$.
Applying $D^c_\lambda$ to $[L_0\ {}_\mu \ Y_i]=(\partial+\mu+\frac{b}{2})Y_i$ and equating coefficients, we obtain
\begin{eqnarray}
\nonumber
\mbox{$(\partial+\lambda+\mu+\frac{b}{2})h_{1,i}(\partial,\lambda)$}&\!\!=\!\!&(\partial+2\mu)h_{1,i}(\partial+\mu,\lambda),\\
\nonumber
\mbox{$(\partial+\lambda+\mu+\frac{b}{2})h_{2,i}(\partial,\lambda)$}&\!\!=\!\!&(\partial+b)h_{2,i}(\partial+\mu,\lambda)+p(\partial+2\lambda+2\mu),\\
\nonumber
\mbox{$(\partial+\lambda+\mu+\frac{b}{2})h_{3,i}(\partial,\lambda)$}&\!\!=\!\!&\mbox{$(\partial+\mu+\frac{b}{2})h_{3,i}(\partial+\mu,\lambda)$}.
\end{eqnarray}
Putting $\mu=0$ in the above equalities, one can deduce that
$h_{1,i}(\partial,\lambda)=h_{2,i}(\partial,\lambda)=h_{3,i}(\partial,\lambda)=p=0$.
Therefore, $D^c_\lambda(L_i)=D^c_\lambda(Y_i)=0$.
Similarly, one can deduce $D^c_\lambda(M_i)=0$ by $[L_0\ {}_\mu \ M_i]=(\partial+b)M_i$.

\vskip2pt
{\bf Case 3:} $a\neq 0, 1$.
\vskip2pt

By \eqref{405} and \eqref{406}, we can set $h_2(\lambda)=\frac{f_{2,0}(-\lambda,\lambda)}{\lambda-b}$ and $h_3(\lambda)=\frac{f_{3,0}(-\lambda,\lambda)}{\lambda-\frac{b}{2}}$.
Replacing $D^c$ by $D^c-{\rm ad}_{h_2(-\partial)M_c+h_3(-\partial)Y_c}$, we may assume that $D^c_\lambda(L_i)=0$.
As above, one can deduce $D^c_\lambda(M_i)=D^c_\lambda(Y_i)=0$ by
relations $[L_0\ {}_\mu \ M_i]=(\partial+a\mu+b)M_i$ and $[L_0\ {}_\mu \ Y_i]=(\partial+(\frac{a}{2}+1)\mu+\frac{b}{2})Y_i$.
\QED\vskip5pt

Now, we can give the proof of Theorem~\ref{main-results-1}.
\vskip5pt

\ni{\it Proof of Theorem \ref{main-results-1}.}\ \
For $D\in {\rm CDer\,} (CSV(a,b))$ and $i\in\Z$, define $D^i$ by
$$
D^i(L_j)=\pi_{i+j} D(L_j),\ \ \ D^i(M_j)=\pi_{i+j} D(M_j),\ \ \ D^i(Y_j)=\pi_{i+j} D(Y_j)\ \ \ \text{for}\ \ \  j\in\Z.
$$
Here, in general, $\pi_{i}$ denotes the natural projection from
$$
\C[\lambda]\otimes CSV(a,b)\cong \big(\OP{k\in\Z}\C[\partial,\lambda]L_k\big)\oplus\big(\OP{k\in\Z}\C[\partial,\lambda]M_k\big)\oplus\big(\OP{k\in\Z}\C[\partial,\lambda]Y_k\big),
$$
onto $\C[\partial,\lambda]{L_{i}}\oplus\C[\partial,\lambda]{M_{i}}\oplus\C[\partial,\lambda]{Y_{i}}$.
Then $D^i$ is a conformal derivation and $D=\mbox{$\sum_{i\in\Z}$} D^i$ in the sense that for any $x\in CSV(a,b)$ only finitely many $D^i_\lambda(x)\neq0$.

By Lemma~\ref{graded-derivations}, we have $D=\sum_{c\in\Z} D^c$, where $D^c=\mbox{ad}_{f(\partial)L_c+g(\partial)M_c+h(\partial)Y_c}+\delta_{a,1}D_{\vec{a}}$ for some
$f(\partial), g(\partial), h(\partial) \in\C[\partial]$ and $\vec{a}=(\ldots, 0, \ldots, 0, q_c, 0, \ldots, 0, \ldots)\in\C^\infty$ with $q_c\in\C$ appearing in the $c$-th position.
Note that
$$
\mbox{$D^c{}_\lambda(L_0)=(\partial+2\lambda)f(-\lambda)L_c
+((a-1)\partial+a\lambda-b)g(-\lambda)M_c
+(\frac{a}{2}\partial+(\frac{a}{2}+1)\lambda-\frac{b}{2})h(-\lambda)Y_c+\delta_{a,1}q_c M_c.$}
$$
It is not difficult to see that $D^c{}_\lambda(L_0)=0$ implies $D^c=0$.
Hence, $D=\sum_{c\in\Z} D^c$ is a finite sum. This completes the proof.
\QED

\section{Rank one conformal modules over $CSV(a,b)$}

In this section, we aim to classify the free conformal modules of rank one over $CSV(a,b)$ for all $a,b\in\C$.
Let us first construct some free conformal modules of rank one over $CSV(a,b)$ for $a,b\in\C$.
Recall that \cite{WCY} a free conformal module of rank one over the loop
Virasoro Lie conformal algebra $CW={\oplus}_{i\in\Z}\C[\partial]L_i$  is isomorphic to $M_{\a,\b,c}$ for some $\a,\b,c\in\C$, where
$M_{\a,\b,c}=\C[\partial]v$ and the $\lambda$-actions are given by
\begin{equation}\label{L-actions}
L_i\,{}_\lambda\, v=c^i(\partial+\a\lambda+\b)v.
\end{equation}
Furthermore, nontrivial $M_{\a,\b,c}$ ($c\ne 0$) is irreducible if and only if $\a\ne 0$.
Obviously, $M_{\a,\b,c}$ is also a free conformal $CSV(a,b)$-module of rank one (still denoted by $M_{\a,\b,c}$)
by extending the $\lambda$-actions of $M_i$ and $Y_i$ trivially, namely
\begin{equation}\label{trivial-extention}
M_i\,{}_\lambda\, v=Y_i\,{}_\lambda\, v=0.
\end{equation}
We focus our interest on whether there are nontrivial extensions of the conformal $CW$-module $M_{\a,\b,c}$. Let us consider $CSV(0,0)$.
For any $d\in\C$, by replacing the $\lambda$-actions \eqref{trivial-extention} by
\begin{equation}\label{nontrivial-extention}
M_i\,{}_\lambda\, v=0,\ \ \ \ Y_i\,{}_\lambda\, v=dc^i v.
\end{equation}
We obtain a free conformal $CSV(0,0)$-module of rank one, denoted by $M_{\a,\b,c,d}$.
Clearly, $M_{\a,\b,c,d}$ with $c, d\ne 0$ is a nontrivial extension of $M_{\a,\b,c}$.
The main result in this section is as follows.

\begin{theo}\label{main-results-2}
Let $M$ be a free conformal module of rank one over $CSV(a,b)$.
\baselineskip1pt\lineskip7pt\parskip-1pt
\begin{itemize}\parskip-1pt
  \item[{\rm(1)}] If $a\ne 0$ or $b\ne 0$, then $M\cong M_{\a,\b,c}$ defined by \eqref{L-actions} and \eqref{trivial-extention} for some $\a,\b,c\in\C$.
  \item[{\rm(2)}] If $a=b=0$, then $M\cong M_{\a,\b,c,d}$ defined by \eqref{L-actions} and \eqref{nontrivial-extention} for some $\a,\b, c, d\in\C$.
\end{itemize}
Furthermore, $M_{\a,\b,c}$ and $M_{\a,\b,c,d}$ are irreducible if and only if $\a, c\ne 0$.
\end{theo}

\ni{\it Proof.}\ \
Write $M=\C[\partial]v$. 
Assume that
$L_i\,{}_\lambda\, v=f_i(\partial,\lambda)v$, $M_i\,{}_\lambda\, v=g_i(\partial,\lambda)v$, $Y_i\,{}_\lambda\, v=h_i(\partial,\lambda)v$, where $f_i(\partial,\lambda)$, $g_i(\partial,\lambda)$, $h_i(\partial,\lambda)\in\C[\partial,\lambda]$.
We only need to determine the coefficients $f_i(\partial,\lambda),g_i(\partial,\lambda),h_i(\partial,\lambda)$.

First, regarding $M$ as a conformal module over $CW$, by \cite{WCY}, we know that $f_i(\partial,\lambda)$
must be of the following form:
\begin{eqnarray}\label{f-form}
f_i(\partial,\lambda)=c^i(\partial+\a\lambda+\b),
\end{eqnarray}
where $\a,\b, c\in\C$.
By Definition~\ref{def-module}, we have
\begin{eqnarray}
\nonumber
M_i\,{}_\lambda\,(M_j\, {}_\mu v)-M_j\,{}_\mu\, (M_i \,{}_\lambda v)&\!\!=\!\!&[M_i\,{}_\lambda\, M_j]\,{}_{\lambda+\mu} v=0,\\
\nonumber
M_i\,{}_\lambda\,(Y_j\, {}_\mu v)-Y_j\,{}_\mu\, (M_i \,{}_\lambda v)&\!\!=\!\!&[M_i\,{}_\lambda\, Y_j]\,{}_{\lambda+\mu} v=0,\\
\nonumber
Y_i\,{}_\lambda\, (Y_j\, {}_\mu v)-Y_j\,{}_\mu\, (Y_i\, {}_\lambda  v)&\!\!=\!\!&[Y_i\,{}_\lambda\, Y_j]\,{}_{\lambda+\mu} v=((\partial+2\lambda)M_{i+j})_{\lambda+\mu}v,
\end{eqnarray}
from which we can get respectively
\begin{eqnarray}
\label{605}
g_j(\partial+\lambda,\mu)g_i(\partial,\lambda)&\!\!=\!\!&g_i(\partial+\mu,\lambda)g_j(\partial,\mu),\\
\label{606}
h_j(\partial+\lambda,\mu)g_i(\partial,\lambda)&\!\!=\!\!&g_i(\partial+\mu,\lambda)h_j(\partial,\mu),\\
\label{607}
h_j(\partial+\lambda,\mu)h_i(\partial,\lambda)-h_i(\partial+\mu,\lambda)h_j(\partial,\mu)&\!\!=\!\!&(\lambda-\mu)g_{i+j}(\partial,\lambda+\mu).
\end{eqnarray}
By comparing the coefficients of $\lambda$ in \eqref{605}, we see that $g_i(\partial,\lambda)$ is independent of the variable $\partial$,
and so we can denote $g_i(\lambda)=g_i(\partial,\lambda)$ for $i\in\Z$.
Then, by \eqref{606}, we see that $h_i(\partial,\lambda)$ is also independent of the variable $\partial$,
and so we can denote $h_i(\lambda)=h_i(\partial,\lambda)$ for $i\in\Z$. Furthermore, by \eqref{607}, we must have
\begin{eqnarray}\label{g-form}
g_i(\partial,\lambda)=0.
\end{eqnarray}
At last, let us further determine the coefficients $h_i(\partial,\lambda)=h_i(\lambda)$.
Since $L_i\,{}_\lambda\,(Y_j\, {}_\mu v)-Y_j\,{}_\mu\, (L_i \,{}_\lambda v)=[L_i\,{}_\lambda\, Y_j]\,{}_{\lambda+\mu} v$, a direct computation shows that
\begin{eqnarray}\label{609}
\mbox{$(\frac{a}{2}\lambda-\mu+\frac{b}{2}) h_{i+j}(\lambda+\mu)=-\mu c^i h_{j}(\mu).$}
\end{eqnarray}

If $a\neq0$ or $b\neq0$, putting $\mu=0$ in \eqref{609}, one can easily deduce that $h_i(\partial,\lambda)=h_i(\lambda)=0$.
This, together with \eqref{f-form} and \eqref{g-form}, implies that the conclusion (1) is true.

If $a=b=0$, putting $i=0$ in \eqref{609}, one can deduce that $h_{j}(\lambda+\mu)= h_{j}(\mu)$. Hence, $h_i(\lambda)$ is independent of the variable $\lambda$, and so
we can denote $h_i=h_i(\lambda)$ for $i\in\Z$.
Then, putting $j=0$ in \eqref{609} with $a=b=0$ and $\mu=1$, we obtain $h_{i}= c^i h_{0}$. Setting $d=h_0$, we have $Y_i\, {}_\lambda \,v=h_i v=d c^i v$.
This, together with \eqref{f-form} and \eqref{g-form}, implies that the conclusion (2) is true.

Clearly, $M_{\a,\b,c}$ or $M_{\a,\b,c,d}$ is nontrivial if and only if $c\ne 0$. Furthermore, as in \cite{WCY} for conformal $CW$-modules,
the nontrivial conformal module $M_{\a,\b,c}$ or $M_{\a,\b,c,d}$ is irreducible if and only if $\a\ne 0$.
\QED

\section{$\Z$-graded free intermediate series modules over $CSV(a,b)$}
In this section, we classify the $\Z$-graded free intermediate series modules over $CSV(a,b)$ for all $a,b\in\C$.
First, we construct some $\Z$-graded free intermediate series modules over $CSV(a,b)$ for $a,b\in\C$.
From \cite{WCY}, we known that there are, up to isomorphism, two classes of nontrivial $\Z$-graded free intermediate series modules over the loop Virasoro Lie conformal algebra
$CW={\oplus}_{i\in\Z}\C[\partial]L_i$. The first class of such modules is $V_{\a,\b}=\oplus_{m\in\Z}\C[\partial]v_m$ ($\a,\b\in\C$) and the $\lambda$-actions are given by
\begin{equation}\label{CW-first}
L_i\,{}_\lambda\, v_m=(\partial+\a\lambda+\b)v_{i+m}.
\end{equation}
Denote by $\{0,1\}^\infty$ the set of sequences
$A=\{a_i\}_{i\in\Z}$  with $a_i\in\{0,1\}$ for any $i\in\Z$. The second class of such modules is
$V_{A,\b}=\oplus_{m\in\Z}\C[\partial]v_m$ ($A\in\{0,1\}^\infty$, $\b\in\C$) and the $\lambda$-actions are given by
\begin{equation}\label{CW-second}
L_i\,{}_{{}_\lambda} v_m=
\begin{cases}
(\partial+\b)v_{i+m} &\ \mbox{if}\  \  \ (a_m,a_{i+m})=(0,0),\\[4pt]
(\partial+\b+\lambda)v_{i+m}&\  \mbox{if} \ \ \ (a_m,a_{i+m})=(1,1),\\[4pt]
v_{i+m} &\ \mbox{if} \ \ \ (a_m,a_{i+m})=(0,1),\\[4pt]
(\partial+\b)(\partial+\b+\lambda)v_{i+m}&\  \mbox{if} \ \ \ (a_m,a_{i+m})=(1,0).
\end{cases}
\end{equation}
Obviously, $V_{\a,\b}$ or $V_{A,\b}$ is also a nontrivial $\Z$-graded free intermediate series module over $CSV(a,b)$ (still denoted by  $V_{\a,\b}$ or $V_{A,\b}$)
by extending the $\lambda$-actions of $M_i$ and $Y_i$ trivially, namely
\begin{equation}\label{trivial-extention-new}
M_i\,{}_\lambda\, v_m=Y_i\,{}_\lambda\, v_m=0.
\end{equation}
As in Section~4, we focus our interest on whether there are nontrivial extensions of the conformal $CW$-modules $V_{\a,\b}$ or $V_{A,\b}$.
Let us consider $CSV(0,0)$. For any $d\in\C$, by replacing the $\lambda$-actions \eqref{trivial-extention-new} by
\begin{equation}\label{nontrivial-extention-new}
M_i\,{}_\lambda\, v_m=0,\ \ \ \ Y_i\,{}_\lambda\, v_m=d v_{i+m}.
\end{equation}
We obtain a nontrivial $\Z$-graded free intermediate series module over $CSV(0,0)$, denoted by $V_{\a,\b,d}$ or $V_{A,\b,d}$.
Clearly, $V_{\a,\b,d}$ (resp.~$V_{A,\b,d}$) with $d\ne 0$ is a nontrivial extension of $V_{\a,\b}$ (resp.~$V_{A,\b}$).
The main result in this section is as follows.

\begin{theo}\label{main-results-3}
Let $V$ be a nontrivial $\Z$-graded free intermediate series module over $CSV(a,b)$.
\baselineskip1pt\lineskip7pt\parskip-1pt
\begin{itemize}\parskip-1pt
  \item[{\rm(1)}] If $a\ne 0$ or $b\ne 0$, then $V\cong V_{\a,\b}$ defined by \eqref{CW-first} and \eqref{trivial-extention-new},
  or $V_{A,\b}$ defined by \eqref{CW-second} and \eqref{trivial-extention-new} for some $A\in\{0,1\}^\infty,\,\a,\b\in\C$.
  \item[{\rm(2)}] If $a=b=0$, then $V\cong V_{\a,\b,d}$ defined by \eqref{CW-first} and \eqref{nontrivial-extention-new},
  or $V_{A,\b,d}$ defined by \eqref{CW-second} and \eqref{nontrivial-extention-new} for some $A\in\{0,1\}^\infty,\,\a,\b,d\in\C$.
\end{itemize}
\end{theo}

Let $V$ be a nontrivial $\Z$-graded free intermediate series module over $CSV(a,b)$.
Write $V=\oplus_{m\in\Z}\C[\partial]v_m$. For $i,j,k,m\in\Z$, we assume that
\begin{eqnarray*}\label{700}
L_i\,{}_\lambda\, v_m=f_{i,m}(\partial,\lambda)v_{i+m}, \ \  M_j\,{}_\lambda\, v_m=g_{j,m}(\partial,\lambda)v_{j+m}, \ \ Y_k\,{}_\lambda\, v_m=h_{k,m}(\partial,\lambda)v_{k+m},
\end{eqnarray*}
where $f_{i,m}(\partial,\lambda)$, $g_{j,m}(\partial,\lambda)$, $h_{k,m}(\partial,\lambda)\in\C[\partial,\lambda]$ are called the {\it structure
coefficients} of $V$  related to the $\C[\partial]$-basis $\{v_m\,|\,m\in\Z\}$.
To prove Theorem~\ref{main-results-3}, we only need to specify these structure
coefficients of $V$.

\begin{lemm}\label{lemm-1}
The structure coefficients of $V$ satisfy the following relations:
\begin{eqnarray}
\label{LM}
g_{j,m}(\partial+\lambda,\mu)f_{i,j+m}(\partial,\lambda)-f_{i,m}(\partial+\mu,\lambda)g_{j,i+m}(\partial,\mu)&\!\!=\!\!&((a-1)\lambda-\mu+b) g_{i+j,m}(\partial,\lambda+\mu),\\
\label{LY}
h_{j,m}(\partial+\lambda,\mu)f_{i,j+m}(\partial,\lambda)-f_{i,m}(\partial+\mu,\lambda)h_{j,i+m}(\partial,\mu)&\!\!=\!\!&\mbox{$(\frac{a}{2}\lambda-\mu+\frac{b}{2})h_{i+j,m}(\partial,\lambda+\mu)$},\\
\label{YY}
h_{j,m}(\partial+\lambda,\mu)h_{i,j+m}(\partial,\lambda)-h_{i,m}(\partial+\mu,\lambda)h_{j,i+m}(\partial,\mu)&\!\!=\!\!&(\lambda-\mu)g_{i+j,m}(\partial,\lambda+\mu),\\
\label{MY}
h_{j,m}(\partial+\lambda,\mu)g_{i,j+m}(\partial,\lambda)&\!\!=\!\!&g_{i,m}(\partial+\mu,\lambda)h_{j,i+m}(\partial,\mu),\\
\label{MM}
g_{j,m}(\partial+\lambda,\mu)g_{i,j+m}(\partial,\lambda)&\!\!=\!\!&g_{i,m}(\partial+\mu,\lambda)g_{j,i+m}(\partial,\mu).
\end{eqnarray}
\end{lemm}
\ni\ni{\it Proof.}\ \
By Definition~\ref{def-module}, we have
\begin{eqnarray*}
[L_i\,{}_\lambda\, M_j]\,{}_{\lambda+\mu} \,v_m &\!\!=\!\!& L_i\,{}_\lambda \,(M_j\, {}_\mu \,v_m)-M_j\,{}_\mu \,(L_i \,{}_\lambda\, v_m) \\
 &\!\!=\!\!& g_{j,m}(\partial+\lambda,\mu)f_{i,j+m}(\partial,\lambda)-f_{i,m}(\partial+\mu,\lambda)g_{j,i+m}(\partial,\mu),
\end{eqnarray*}
and
\begin{equation*}\label{709}
((\partial+a\lambda+b)M_{i+j})\,{}_{\lambda+\mu} \,v_m=((a-1)\lambda-\mu+b) g_{i+j,m}(\partial,\lambda+\mu).
\end{equation*}
Then the relation \eqref{LM} follows from \eqref{brackets-2}.
Similarly, the relations \eqref{LY}--\eqref{MM} can be deduced from  \eqref{brackets-3}, \eqref{brackets-4} and
$[M_i\,{}_\lambda\, Y_j]=[M_i\,{}_\lambda\, M_j]=0$.
\QED

\begin{lemm}\label{lemm-2}
Suppose $V\cong V_{\a,\b}$ as $CW$-modules for some $\a,\b\in\C$.
\baselineskip1pt\lineskip7pt\parskip-1pt
\begin{itemize}\parskip-1pt
  \item[{\rm(1)}] If $a\ne 0$ or $b\ne 0$, then $g_{i,m}(\partial,\lambda)=h_{i,m}(\partial,\lambda)=0$.
  \item[{\rm(2)}] If $a=b=0$, then $g_{i,m}(\partial,\lambda)=0$ and $h_{i,m}(\partial,\lambda)=d$ for some $d\in\C$.
\end{itemize}
\end{lemm}
\ni{\it Proof.}\ \
First, by \cite{WCY}, we have $f_{i,m}(\partial,\lambda)=\partial+\a\lambda+\b$.
As in the proof of Theorem~\ref{main-results-2}, from \eqref{MM} with $i=j=0$, we see that $g_{0,m}(\partial,\lambda)$ is independent of the variable $\partial$.
Thus, we can denote $g_{0,m}(\lambda)=g_{0,m}(\partial,\lambda)$.
Taking $i=j=0$ in \eqref{LM}, we obtain
\begin{equation}\label{714}
-\mu g_{0,m}(\mu)=((a-1)\lambda-\mu+b) g_{0,m}(\lambda+\mu).
\end{equation}
We divide the remaining proof into the following four cases.

\vskip2pt
{\bf Case 1:} $b\neq0$.
\vskip2pt

In this case, taking $\lambda=0$ in \eqref{714} and noting that $b\ne0$, we see that $g_{0,m}(\partial,\lambda)=g_{0,m}(\lambda)=0$.
Then, by \eqref{LM} with $j=0$, we obtain $((a-1)\lambda-\mu+b) g_{i,m}(\partial,\lambda+\mu)=0$, which implies
$g_{i,m}(\partial,\lambda)=0$.

Now, \eqref{YY} becomes
\begin{equation}\label{717-new}
h_{j,m}(\partial+\lambda,\mu)h_{i,j+m}(\partial,\lambda)=h_{i,m}(\partial+\mu,\lambda)h_{j,i+m}(\partial,\mu).
\end{equation}
Taking $i=j=0$ in \eqref{717-new} and comparing the coefficients of $\lambda$, we see that $h_{0,m}(\partial,\lambda)$ is independent of the variable $\partial$, and so
we can denote $h_{0,m}(\lambda)=h_{0,m}(\partial,\lambda)$.
Then, taking $i=j=0$ in \eqref{LY}, we obtain
\begin{equation}\label{719}
\mbox{$(\frac{a}{2}\lambda-\mu+\frac{b}{2})h_{0,m}(\lambda+\mu)=-\mu h_{0,m}(\mu).$}
\end{equation}
Furthermore, taking $\lambda=0$ in \eqref{719} and noting that $b\ne0$, we see that $h_{0,m}(\partial,\lambda)=h_{0,m}(\lambda)=0$.
Then, by \eqref{LY} with $j=0$, one can deduce that $h_{i,m}(\partial,\lambda)=0$.

\vskip2pt
{\bf Case 2:} $a=b=0$.
\vskip2pt

In this case, by taking $\mu=0$ in \eqref{714}, we see that $g_{0,m}(\partial,\lambda)=g_{0,m}(\lambda)=0$.
Then, by \eqref{LM} with $j=0$, we obtain $(\lambda+\mu) g_{i,m}(\partial,\lambda+\mu)=0$, which implies
$g_{i,m}(\partial,\lambda)=0$.

As in Case~1, by \eqref{YY}, one can deduce that $h_{0,m}(\partial,\lambda)$ is independent of the variable $\partial$.
Thus, we can denote $h_{0,m}(\lambda)=h_{0,m}(\partial,\lambda)$.
Then, taking $i=j=0$ in \eqref{LY}, we obtain
$h_{0,m}(\lambda+\mu)=h_{0,m}(\mu)$, which implies that $h_{0,m}(\lambda)$ is independent of the variable $\lambda$, and so
we can further denote $h_m=h_{0,m}(\lambda)=h_{0,m}(\partial,\lambda)$.
Again, taking $j=0$ in \eqref{LY},  we obtain
\begin{equation}\label{723}
\mu h_{i,m}(\partial,\lambda+\mu)=(\partial+\a\lambda+\b)(h_{i+m}-h_{m})+\mu h_{i+m}.
\end{equation}
Equating the coefficients of $\mu^n$ for $n=0,1$ in \eqref{723},
we immediately obtain $h_m=d$ for some $d\in\C$, and $h_{i,m}(\partial,\lambda)=d$.

\vskip2pt
{\bf Case 3:} $a=1$, $b=0$.
\vskip2pt

In this case, by \eqref{714}, we have $g_{0,m}(\mu)=g_{0,m}(\lambda+\mu)$. Thus, $g_{0,m}(\lambda)$ is independent of the variable $\lambda$, and so
we can denote $g_m=g_{0,m}(\lambda)=g_{0,m}(\partial,\lambda)$. Then, taking $j=0$ in \eqref{LM}, we obtain
\begin{equation}\label{723---new}
\mu g_{i,m}(\partial,\lambda+\mu)=(\partial+\a\lambda+\b)(g_{i+m}-g_{m})+\mu g_{i+m}.
\end{equation}
As in Case~2, by equating the coefficients of $\mu^n$ for $n=0,1$ in \eqref{723---new},
one can deduce that there exists some $e\in\C$ such that $g_{i,m}(\partial,\lambda)=e$.
We claim that $e=0$. In fact, if this is not true, by \eqref{MY} we see that
$h_{j,m}(\partial+\lambda,\mu)=h_{j,i+m}(\partial,\mu)$, which implies that
$h_{i,m}(\partial,\lambda)$ is independent of the variable $\partial$ and the subscript $m$. Then, by \eqref{YY}, we obtain
$0=(\lambda-\mu)e$, a contradiction. Hence, $g_{i,m}(\partial,\lambda)=0$.

Now, as in Case~1, by \eqref{YY}, one can deduce that $h_{0,m}(\partial,\lambda)$ is independent of the variable $\partial$.
Thus, we can denote $h_{0,m}(\lambda)=h_{0,m}(\partial,\lambda)$.
Then, taking $i=j=0$ in \eqref{LY}, we obtain
\begin{equation}\label{719-----new}
\mbox{$(\frac{1}{2}\lambda-\mu)h_{0,m}(\lambda+\mu)=-\mu h_{0,m}(\mu).$}
\end{equation}
Furthermore, taking $\mu=0$ in \eqref{719-----new}, we see that $h_{0,m}(\partial,\lambda)=h_{0,m}(\lambda)=0$.
Then, by \eqref{LY} with $j=0$, one can deduce that $h_{i,m}(\partial,\lambda)=0$.

One can also prove the results in this case base on \cite[Theorem~6.6]{FSW} (see \cite{CHSX}).

\vskip2pt
{\bf Case 4:} $a\neq 0,1$, $b=0$,
\vskip2pt

 As in Case~2, taking $\mu=0$ in \eqref{714} and noting that $a\ne 1$, we see that $g_{0,m}(\partial,\lambda)=g_{0,m}(\lambda)=0$.
Then, by \eqref{LM} with $j=0$, we obtain $((a-1)\lambda-\mu) g_{i,m}(\partial,\lambda+\mu)=0$, which implies
$g_{i,m}(\partial,\lambda)=0$.

As in Case~1, by \eqref{YY}, one can deduce that $h_{0,m}(\partial,\lambda)$ is independent of the variable $\partial$.
Thus, we can denote $h_{0,m}(\lambda)=h_{0,m}(\partial,\lambda)$. Then, taking $i=j=0$ in \eqref{LY}, we obtain
\begin{equation}\label{719--------new}
\mbox{$(\frac{a}{2}\lambda-\mu)h_{0,m}(\lambda+\mu)=-\mu h_{0,m}(\mu).$}
\end{equation}
Furthermore, taking $\mu=0$ in \eqref{719--------new} and noting that $a\ne0$, we see that $h_{0,m}(\partial,\lambda)=h_{0,m}(\lambda)=0$.
Then, by \eqref{LY} with $j=0$, one can deduce that $h_{i,m}(\partial,\lambda)=0$.
\QED

\vs{15pt}

Similar to the above lemma, we can also prove the following lemma.

\begin{lemm}\label{lemm-3}
Suppose $V\cong V_{A,\b, d}$ as $CW$-modules for some $A\in\{0,1\}^\infty$, $\b,d\in\C$.
\baselineskip1pt\lineskip7pt\parskip-1pt
\begin{itemize}\parskip-1pt
  \item[{\rm(1)}] If $a\ne 0$ or $b\ne 0$, then $g_{i,m}(\partial,\lambda)=h_{i,m}(\partial,\lambda)=0$.
  \item[{\rm(2)}] If $a=b=0$, then $g_{i,m}(\partial,\lambda)=0$ and $h_{i,m}(\partial,\lambda)=d$ for some $d\in\C$.
\end{itemize}
\end{lemm}

Now, Theorem~\ref{main-results-3} follows from Lemmas~\ref{lemm-2} and \ref{lemm-3}.

\section{Results of Lie conformal algebras $CHV(a,b)$}

In this section, we present the results of the infinite rank Heisenberg-Virasoro type Lie conformal algebra $CHV(a,b)=({\oplus}_{i\in\Z}\C[\partial]L_i)\oplus({\oplus}_{i\in\Z}\C[\partial]M_i)$ for $a,b\in\C$.
The proofs are similar to those for $CSV(a,b)$ in previous sections; the details are omitted.

(I) If $a\ne 1$, then all conformal derivations of $CHV(a,b)$ are inner. If $a=1$, then there exist non-inner conformal derivations of $CHV(1,b)$.
Suppose that $D$ is a non-inner conformal derivation of $CHV(1,b)$. Then $D=D_{\vec{a}}$ for some (non-zero) $\vec{a}\in \C^\infty$, where $D_{\vec{a}}$ is defined by
\begin{equation*}
D_{\vec{a}}\, {}_\lambda\, (L_i)=\sum_{c\in\Z} a_c M_{i+c}, \ \ \ D_{\vec{a}}\, {}_\lambda\,(M_i)=0 \ \ \ \mbox{for}\ \ \  i\in\Z.
\end{equation*}
Denote by ${\rm CDer\,}(CHV(a,b))$ and ${\rm CInn\,}(CHV(a,b))$ the vector spaces of all conformal derivations and inner conformal derivations of $CHV(a,b)$, respectively.

\begin{theo}\label{main-results-1-CHV}
We have ${\rm CDer\,}(CHV(a,b))={\rm CInn\,}(CHV(a,b))\oplus \delta_{a,1}\C^\infty$.
\end{theo}

(II) If $a\neq 1$ or $b\neq 0$, then a free conformal $CHV(a,b)$-module of rank one
is simply a free conformal $CW$-module of rank one $M_{\a,\b,c}=\C[\partial]v$ (c.f.~\eqref{L-actions})
with trivial $\lambda$-actions of $M_i$, namely
\begin{equation}\label{trivial-extention-CHV}
M_i\,{}_\lambda\, v=0.
\end{equation}
Denote such a module by $M'_{\a,\b,c}$.
If $a=1,\,b=0$, then, for any $d\in\C$, replacing the $\lambda$-actions \eqref{trivial-extention-CHV} by
\begin{equation}\label{nontrivial-extention-CHV}
M_i\,{}_\lambda\, v=dc^i v,
\end{equation}
we obtain a free conformal $CHV(1,0)$-module of rank one, denoted by $M'_{\a,\b,c,d}$.
Note that $M'_{\a,\b,c,d}$ is a nontrivial extension of $M_{\a,\b,c}$ if and only if $c, d\ne 0$.

\begin{theo}\label{main-results-2-CHV}
Let $M$ be a free conformal module of rank one over $CHV(a,b)$.
\baselineskip1pt\lineskip7pt\parskip-1pt
\begin{itemize}\parskip-1pt
  \item[{\rm(1)}] If $a\ne 1$ or $b\ne 0$, then $M\cong M'_{\a,\b,c}$ defined by \eqref{L-actions} and \eqref{trivial-extention-CHV} for some $\a,\b,c\in\C$.
  \item[{\rm(2)}] If $a=1,\,b=0$, then $M\cong M'_{\a,\b,c,d}$ defined by \eqref{L-actions} and \eqref{nontrivial-extention-CHV}  for some $\a,\b, c, d\in\C$.
\end{itemize}
Furthermore, $M'_{\a,\b,c}$ and $M'_{\a,\b,c,d}$ are irreducible if and only if $\a, c\ne 0$.
\end{theo}

(III) If $a\neq 1$ or $b\neq 0$, then a nontrivial $\Z$-graded free intermediate series module over $CHV(a,b)$
is simply a nontrivial $\Z$-graded free intermediate series $CW$-module $V_{\a,\b}=\oplus_{m\in\Z}\C[\partial]v_m$ (c.f.~\eqref{CW-first})
or $V_{A,\b}=\oplus_{m\in\Z}\C[\partial]v_m$ (c.f.~\eqref{CW-second}) with trivial $\lambda$-actions of $M_i$, namely
\begin{equation}\label{trivial-extention-new-CHV}
M_i\,{}_\lambda\, v_m=0.
\end{equation}
Denote such a module by $V'_{\a,\b}$ or $V'_{A,\b}$.
If $a=1,\,b=0$, then for any $d\in\C$, replacing the $\lambda$-actions \eqref{trivial-extention-new-CHV} by
\begin{equation}\label{nontrivial-extention-new-CHV}
M_i\,{}_\lambda\, v_m=d v_{i+m},
\end{equation}
we obtain a nontrivial $\Z$-graded free intermediate series $CHV(1,0)$-module, denoted by $V'_{\a,\b,d}$ or $V'_{A,\b,d}$.
Note that $V'_{\a,\b,d}$ (resp.~$V'_{A,\b,d}$) is a nontrivial extension of $V_{\a,\b}$ (resp.~$V_{A,\b}$) if and only if $d\ne 0$.

\begin{theo}\label{main-results-3-CHV}
Let $V$ be a nontrivial $\Z$-graded free intermediate series module over $CHV(a,b)$.
\baselineskip1pt\lineskip7pt\parskip-1pt
\begin{itemize}\parskip-1pt
  \item[{\rm(1)}] If $a\ne 1$ or $b\ne 0$, then $V\cong V'_{\a,\b}$ defined by \eqref{CW-first} and \eqref{trivial-extention-new-CHV},
  or $V'_{A,\b}$ defined by \eqref{CW-second} and \eqref{trivial-extention-new-CHV} for some $A\in\{0,1\}^\infty,\,\a,\b\in\C$.
  \item[{\rm(2)}] If $a=1,\,b=0$, then $V\cong V'_{\a,\b,d}$ defined by \eqref{CW-first} and \eqref{nontrivial-extention-new-CHV},
  or $V'_{A,\b,d}$ defined by \eqref{CW-second} and \eqref{nontrivial-extention-new-CHV} for some $A\in\{0,1\}^\infty,\,\a,\b,d\in\C$.
\end{itemize}
\end{theo}

\vskip10pt

\small \ni{\bf Acknowledgement}
This work was supported by the National Natural Science Foundation of China (Nos.~11431010, 11371278, 11401570),
the Natural Science Foundation of Jiangsu Province, China (No.~BK20140177),
and the Fundamental Research Funds for the Central Universities (No.~2014QNA68).

\end{document}